\documentclass[12pt,a4paper] {article}

\usepackage[utf8]{inputenc}
\usepackage[english]{babel}

\usepackage{amssymb}

\usepackage[namelimits,intlimits] {amsmath}
\usepackage{amsthm}

\begin{document}

\title {THE PROCEDURE OF EXCLUDING OF THE NONLINEAR TREND FOR THE MODELS DESCRIBED BY STOCHASTIC DIFFERENTIAL AND DIFFERENCE EQUATIONS\thanks{The article was prepared within the framework of the Academic Fund Program at the National Research University Higher School of Economics (HSE) and supported within the framework of a subsidy granted to the HSE by the Government of the Russian Federation for the implementation of the Global Competitiveness Program}}%
\date{} 
\author {V.~KONAKOV \;\;\; A.~MARKOVA\\National Research University Higher School of Economics (HSE) Moscow}
\maketitle

\begin{abstract}
We consider the diffusion process and its approximation by Markov chain with nonlinear increasing trends. The usual parametrix method is not appliable because these models have unbounded trends. We describe a procedure that allows to exclude nonlinear growing trend and move to stochastic differential equation with reduced drift and diffusion coefficients. A similar procedure is considered for a Markov chain 
\\
\end{abstract}

\section{Introduction}

This article is a continuation of author's article \cite{KMa15}, which considered the case of a linear trend. It is well known that the parametrix method works well for diffusions with bounded drift and diffusion coefficients (see \cite{IKO62}, \cite{Fr64}, \cite{K12}). But in practice we often deal with unbounded coefficients. The aim of the present note is to consider diffusions with unbounded drift and to propose a procedure which reduces the initial equation with unbounded drift to the equation with bounded drift. Similar procedure may be developed for Markov chains. The idea is to consider a new process and then to write a stochastic differential for this process using the Ito formula. To exclude the drift we return back along the solutions of some ODE. This ODE is obtained by discarding the Brownian part and the bounded part of the drift of the initial SDE. To apply the Ito formula we need smoothness w.r.t. initial conditions and initial time. That is why we need corresponding smoothnes assumptions on the unbounded part of the drift \cite{Ar98}. Note also that to apply our approach we need existence and uniqueness results for ODE not only locally but on the whole interval $[0,T]$. Of course, local approach may be developed inder weaker conditions.

In what follows the asterisk like $F_{*}(t,\mathbf{x)}$ denotes the derivative with respect to the space variable $\mathbf{x}$ for fixed $t$. Thus $F_{*}(t,\mathbf{x)}$ is a linear operator from $\mathbb{R}^{d}$ into $\mathbb{R}^{d}$. We use also the standard notation for the higher order derivatives $D_{x}^{\nu }F(t,x)$.

	\section{Assumptions and main results}
	
	We consider the following diffusion model:
	\begin{eqnarray}\label{1}
	dY_{t} &=& \left\{ F(t,Y_{t})+m(t,Y_{t})\right\} dt+\sigma (t,Y_{t})dW_{t},
	\\ \nonumber
	Y_{0}&=&x_{0}\in R^{d}, \;\;\;\; 0\leq t\leq T,  
	\end{eqnarray}
	where $F(t,y)$ and $m(t,y)$ id $d$-dimentional vector-functions, $W_{t}$ -- $d$-dimentional Wiener process. Let us introduce the following assumptions.
	
	\textbf{(A1) (Uniform Ellipticity)}
	The matrix $a:=\sigma \sigma ^{\ast }$ is uniformly elliptic, i.e. there is $\Lambda \geq 1$, $\forall (t,x,\xi )\in
	\lbrack 0,T]\times \left( R^{d}\right) ^{2}$
	\begin{equation*}
	\Lambda ^{-1}\left\vert \xi \right\vert ^{2}\leq \left\langle a(t,x)\xi ,\xi
	\right\rangle \leq \Lambda \left\vert \xi \right\vert ^{2}.
	\end{equation*}
	
	\textbf{(A2) (Smoothness of the drift)}
	The function $F(t,x)$ is two times continuously differentiable in $[0, T]\times R^{d},$ and
	\begin{equation*}
	\max_{t\in \lbrack 0,T]}\left\Vert D_{x}^{\nu }F(t,x)\right\Vert \leq
	M_{F},\left\vert \nu \right\vert =1,2.
	\end{equation*}
	
	\textbf{(A3) (Boundedness of the drift) } The function $m(t,x)$ is continuous and bounded in $(t,x)\in \lbrack 0,T]\times R^{d}$. 
	
	It follows from \textbf{(A2) } and Theorem 1.2 (\cite{DK70}, p.392) that for any  $\mathbf{x}_{0}\in R^{d}$ and $t_{0}\in \lbrack 0,T]$ the ODE
	\begin{equation}
	\overset{\cdot }{\mathbf{x}}=F(t,\mathbf{x})  \label{2}
	\end{equation}
	has unique solution $\mathbf{g}(t):=\mathbf{g(}t;t_{0},\mathbf{x}_{0}\mathbf{)}$ defined on $[0,T]$ such that $\mathbf{g}(t_{0})=\mathbf{x}_{0}$.
	Consider the mapping $G:$ $[0,T]\times R^{d}$ $\rightarrow \lbrack
	0,T]\times R^{d}$ defined by the formula $G(t,\mathbf{x)=(}t,\mathbf{g(}
	t;t_{0},\mathbf{x)})$, where 
	\begin{equation*}
	\mathbf{g(}t;t_{0},\mathbf{x)=}\left( \mathbf{g}_{1}\mathbf{(}t;t_{0},
	\mathbf{x),...,g}_{d}\mathbf{(}t;t_{0},\mathbf{x)}\right) 
	\end{equation*}
	is a solution of the ODE (\ref{2}) satifying the initial condition $\mathbf{g(}
	t_{0};t_{0},\mathbf{x)=x}$. By the refined version of the rectification theorem for nonautonomous case  (\cite{Ar98}, p. 230) and assumption \textbf{(A2)} the mapping $G$ is a local diffeomorphism of class \textbf{\ }$C^{2}$	in a neighborhood of the point $(t_{0},\mathbf{x}_{0}).$ 
	
	The derivative $\mathbf{g}_{\ast}(t;t_{0},\mathbf{x):=}\left\Vert z_{ik}(t;t_{0},\mathbf{x)}
	\right\Vert ,$ $z_{ik}(t;t_{0},\mathbf{x)=}\frac{\partial \mathbf{g}_{i}
		\mathbf{(}t;t_{0},\mathbf{x)}}{\partial x_{k}},i,k=1,...,d,$ of the solution of equation (\ref{2}) with respect to the initial condition $\mathbf{x}$ satisfies the equation of variations with the initial conditon $\mathbf{g}
	_{\ast }(t_{0};t_{0},\mathbf{x)=I}$, where $\mathbf{I}$ is the identity matrix: 
	\begin{equation*}
	\frac{\partial }{\partial t}\mathbf{g(}t;t_{0},\mathbf{x)=}F(t,\mathbf{g(}
	t;t_{0},\mathbf{x)}),  \label{3}
	\end{equation*}
	\begin{equation}
	\frac{\partial }{\partial t}\mathbf{g}_{\ast }(t;t_{0},\mathbf{x)=}F_{\ast
	}(t,\mathbf{g(}t;t_{0},\mathbf{x))\mathbf{g}_{\ast }(}t;t_{0},\mathbf{
	\mathbf{x)},}  \label{4}
\end{equation}
\begin{equation*}
\mathbf{g(}t_{0};t_{0},\mathbf{x)=x,\mathbf{g}_{\ast }(}t_{0};t_{0},\mathbf{
	\mathbf{x)=x}}_{\ast }=\mathbf{I.}  \label{5}
\end{equation*}
(for the details see \cite{Ar98}, p. 225). Moreover,
\begin{equation}
\frac{\partial }{\partial t_{0}}\mathbf{g(}t;t_{0},\mathbf{\mathbf{x})=-
	\mathbf{g}_{\ast }(}t;t_{0},\mathbf{\mathbf{x)}}F(t_{0},\mathbf{\mathbf{x).}}
\label{6}
\end{equation}
For the proof of (\ref{6}) see, e.g. \cite{Bi91}. Clearly, $\mathbf{\mathbf{g}
	_{\ast }(}t;t_{0},\mathbf{\mathbf{x)}}$ is the matrix of the differential at $\mathbf{\mathbf{x}}$ of the $(t_{o},t)$ -- advance transformation given by the phase flow. Analogously, the inverse matrix $\mathbf{\mathbf{g}_{\ast}^{-1}(}t;t_{0},\mathbf{\mathbf{x)}}$ is the matrix of the differential of the $(t,t_{0})$ -- back transformation given by the inverse phase flow.
	Relation (\ref{6}) may be interpreted from geometrical arguments. The derivative of the solution w.r.t. the initial time is the vector image of the tangent vector $-F(t_{0},\mathbf{\mathbf{x)}}$ under the action of the differential at $\mathbf{\mathbf{x}}$ of the $(t_{o},t)$ -- advance transformation given by the phase flow. The matrix of the differential of the inverse mapping is equal to the inverse matrix of the differential of the mapping given by the forward phase flow and symmetrically to (\ref{6}) we obtain:
\begin{equation} \label{7}
\frac{\partial \mathbf{\mathbf{g}}^{-1}(t_{0};t,\mathbf{y\mathbf{)}}}{
	\partial t}=-\mathbf{g}_{\ast }^{-1}(t;t_{0},\mathbf{\mathbf{g}}
^{-1}(t_{0};t,\mathbf{y)})F(t,\mathbf{y}). 
\end{equation}

Note that in the linear case $F(t,\mathbf{x)=}b(t)\mathbf{x}$ we get in notations of \cite{KMa15}:
\begin{equation*}
F(t,\mathbf{g(}t;t_{0},\mathbf{x)})=b(t)\mathbf{g(}t;t_{0},\mathbf{x)},\text{
}F_{\ast }(t,\mathbf{g(}t;t_{0},\mathbf{x))}=b(t),\mathbf{g}_{\ast }(t;t_{0},
\mathbf{x)}=\Phi (t),
\end{equation*}
and from (\ref{4}) and (\ref{7}) we get relations corresponding to the
linear case: 
\begin{equation*}
\Phi ^{\prime }(t)=b(t)\Phi (t),\left[ \Phi ^{-1}(t)\right] ^{\prime }=
\mathbf{-}\left[ \Phi (t)\right] ^{-1}b(t).
\end{equation*}

From now we put $t_{0}=0$. Consider a new stochastic process
\begin{equation*}
\tilde{Y}_{t}=\mathbf{\mathbf{g}}^{-1}(0;t,Y_{t}), \;\;\;\;
\tilde{Y}_{0}=Y_{0}=x_{0},
\end{equation*}
where $Y_{t}$ is a solution of (\ref{1}). Denote $\ \Psi (t,\mathbf{y}):=\mathbf{\mathbf{g}}^{-1}(0;t,\mathbf{y)}$ and apply the $d$--dimensional Ito formula to the process $\widetilde{Y}_{t}=\Psi (t,Y_{t})$. To apply the Ito formula we have to calculate first and second derivatives w.r.t. the terminal values $\mathbf{y}$ and first derivatives w.r.t. the terminal time $t$. Differentiating the identity
\begin{equation*}
\mathbf{g}_{i}\mathbf{(}t;0,\mathbf{\mathbf{g}}^{-1}(0;t,\mathbf{y))=}y_{i},
\end{equation*}
we obtain 
\begin{equation*}
\frac{\partial }{\partial y_{k}}\mathbf{g}_{i}\mathbf{(}t;0,\mathbf{\mathbf{g
	}}^{-1}(0;t,\mathbf{y))}=\sum_{l=1}^{d}\frac{\partial \mathbf{g}_{i}\mathbf{
}t;0,\mathbf{\mathbf{g}}^{-1}(0;t,\mathbf{y))}}{\partial \mathbf{\mathbf{g}}
_{l}^{-1}(0;t,\mathbf{y)}}\frac{\partial \mathbf{\mathbf{g}}_{l}^{-1}(0;t,
\mathbf{y)}}{\partial y_{k}}=\frac{\partial y_{i}}{\partial y_{k}}=\delta
_{ik}  \label{8}
\end{equation*}
and, hence,
\begin{equation}
\mathbf{\mathbf{g}}_{\ast }^{-1}(t;0,\mathbf{\mathbf{g}}^{-1}(0;t,\mathbf{
	y))=}\left\Vert \frac{\partial \mathbf{\mathbf{g}}_{i}^{-1}(0;t,\mathbf{y)}}{
	\partial y_{k}}\right\Vert \mathbf{:=}\left\Vert z^{ik}(t;0,\mathbf{\mathbf{g
	}}^{-1}(0;t,\mathbf{y))}\right\Vert .  \label{9}
	\end{equation}
	From (\ref{4}) and the Liouville's Theorem, see e.g. \cite{Ar98}, for some $C>1$ , we obtain for some
	\begin{equation}
	C^{-1}\leq \det \mathbf{g}_{\ast }(t;t_{0},\mathbf{x)=}\exp \int_{0}^{t}
	\text{trace}[F_{\ast }(s,\mathbf{g(}s;t_{0},\mathbf{x))]}ds\leq C,  \label{10}
	\end{equation}
	and by Theorem 8.65 from \cite{KP04}
	\begin{equation}
	\left\Vert \mathbf{\mathbf{g}_{\ast }(}t;t_{0},\mathbf{\mathbf{x)}}
	\right\Vert \leq C(d,T).  \label{11}
	\end{equation}
	It follows from (\ref{9})-(\ref{11}) that all elements $z^{ik}(t;0,\mathbf{\mathbf{g}}^{-1}(0;t,\mathbf{y}))$ of the inverse matrix $\mathbf{\mathbf{g}}_{\ast }^{-1}(t;0,\mathbf{\mathbf{g}}^{-1}(0;t,
	\mathbf{y))}$ are bounded and because of boundedness of $m(t,Y_{t})$ and $
	\sigma (t,Y_{t})$, the elements of $\widetilde{m}(t,Y_{t}):=\mathbf{\mathbf{g}
	}_{\ast }^{-1}(t;0,\mathbf{\mathbf{g}}^{-1}(0;t,Y_{t}\mathbf{))}m(t,Y_{t})$
	and $\widetilde{\sigma }(t,Y_{t}):=\mathbf{\mathbf{g}}_{\ast }^{-1}(t;0,
	\mathbf{\mathbf{g}}^{-1}(0;t,Y_{t}\mathbf{))}\sigma (t,Y_{t})$ are also bounded. It follows from (\ref{8}) that 
	\begin{equation*}
	\frac{\partial ^{2}}{\partial y_{j}\partial y_{k}}\mathbf{g}_{i}\mathbf{(}
	t;0,\mathbf{\mathbf{g}}^{-1}(0;t,\mathbf{y))}=\sum_{l=1}^{d}\frac{\partial }{
		\partial y_{j}}\left[ \frac{\partial \mathbf{g}_{i}\mathbf{(}t;0,\mathbf{
			\mathbf{g}}^{-1}(0;t,\mathbf{y))}}{\partial \mathbf{\mathbf{g}}_{l}^{-1}(0;t,
		\mathbf{y)}}\right] \frac{\partial \mathbf{\mathbf{g}}_{l}^{-1}(0;t,\mathbf{
			y)}}{\partial y_{k}}+
	\end{equation*}
	\begin{equation*}
	\sum_{l=1}^{d}\frac{\partial \mathbf{g}_{i}\mathbf{(}t;0,\mathbf{\mathbf{g}}
		^{-1}(0;t,\mathbf{y))}}{\partial \mathbf{\mathbf{g}}_{l}^{-1}(0;t,\mathbf{y)}
	}\frac{\partial ^{2}\mathbf{\mathbf{g}}_{l}^{-1}(0;t,\mathbf{y)}}{\partial
	y_{j}\partial y_{k}}=0,
\end{equation*}
or, in our notations,
\begin{equation*}
\frac{\partial ^{2}}{\partial y_{j}\partial y_{k}}\mathbf{g}_{i}\mathbf{(}
t;0,\mathbf{\mathbf{g}}^{-1}(0;t,\mathbf{y))}=\sum_{l,p=1}^{d}\frac{\partial
	z_{il}(t,0,\mathbf{\mathbf{g}}^{-1}(0;t,\mathbf{y)})}{\partial \mathbf{
		\mathbf{g}}_{p}^{-1}(0;t,\mathbf{y)}}z^{pj}(t;0,\mathbf{\mathbf{g}}^{-1}(0;t,
\mathbf{y)})\times 
\end{equation*}
\begin{equation}
z^{lk}(t;0,\mathbf{\mathbf{g}}^{-1}(0;t,\mathbf{y)})+
\sum_{l=1}^{d}z_{il}(t,0,\mathbf{\mathbf{g}}^{-1}(0;t,\mathbf{y)})\frac{
	\partial ^{2}\mathbf{\mathbf{g}}_{l}^{-1}(0;t,\mathbf{y)}}{\partial
	y_{j}\partial y_{k}}=0.  \label{12}
\end{equation}
Differentiating both sides of (\ref{4}) w.r.t. $\mathbf{x}$ and applying \textbf{(A2)} and the Gronwall's Inequality, (see, e.g. \cite{KP04}, Theorem 8.65) we obtain that uniformly in  $(t,\mathbf{y)\in }[0,T]\times R^{d}$ and for $1\leq i,l,p\leq d$ 
\begin{equation*}
\left\vert \frac{\partial z_{il}(t,0,\mathbf{\mathbf{g}}^{-1}(0;t,\mathbf{y)}
	)}{\partial \mathbf{\mathbf{g}}_{p}^{-1}(0;t,\mathbf{y)}}\right\vert \leq C,
\end{equation*}
where the constant $C$ depends only on $T$ and $M_{F}$ from \textbf{(A2).}
Hence 
\begin{equation*}
c_{jk}^{i}(t,0,\mathbf{\mathbf{g}}^{-1}(0;t,\mathbf{y)}):=\sum_{l,p=1}^{d}
\frac{\partial z_{il}(t,0,\mathbf{\mathbf{g}}^{-1}(0;t,\mathbf{y)})}{
	\partial \mathbf{\mathbf{g}}_{p}^{-1}(0;t,\mathbf{y)}}z^{pj}(t,0,\mathbf{
	\mathbf{g}}^{-1}(0;t,\mathbf{y)})\times 
\end{equation*}
\begin{equation*}
z^{lk}(t;0,\mathbf{\mathbf{g}}^{-1}(0;t,\mathbf{y)})\leq C_{1}<\infty .
\end{equation*}
We obtain from (\ref{12}) 
\begin{equation}
\left( 
\begin{array}{c}
\frac{\partial ^{2}\mathbf{\mathbf{g}}_{1}^{-1}(0;t,\mathbf{y)}}{\partial
	y_{j}\partial y_{k}} \\ 
\vdots  \\ 
\frac{\partial ^{2}\mathbf{\mathbf{g}}_{d}^{-1}(0;t,\mathbf{y)}}{\partial
	y_{j}\partial y_{k}}
\end{array}
\right) =-\mathbf{\mathbf{g}}_{\ast }^{-1}(t;0,\mathbf{\mathbf{g}}^{-1}(0;t,
\mathbf{y))}c_{jk}(t,0,\mathbf{\mathbf{g}}^{-1}(0;t,\mathbf{y)}),  \label{13}
\end{equation}
where 
\begin{equation*}
c_{jk}(t,0,\mathbf{\mathbf{g}}^{-1}(0;t,\mathbf{y)})=\left( 
\begin{array}{c}
c_{jk}^{1}(t,0,\mathbf{\mathbf{g}}^{-1}(0;t,\mathbf{y)}) \\ 
\vdots  \\ 
c_{jk}^{d}(t,0,\mathbf{\mathbf{g}}^{-1}(0;t,\mathbf{y)})
\end{array}
\right) ,
\end{equation*}
and,hence,all components of the vector in the l.h.s. of (\ref{13}) are
bounded. By the $d$--dimensional Ito formula:
\begin{equation*}
d\widetilde{Y}_{t,k}=\frac{\partial \mathbf{\mathbf{g}}_{k}^{-1}(0;t,Y_{t}
	\mathbf{)}}{\partial t}dt+\sum_{i=1}^{d}\frac{\partial \mathbf{\mathbf{g}}
	_{k}^{-1}(0;t,Y_{t}\mathbf{)}}{\partial y_{i}}dY_{t,i}+
\end{equation*}
\begin{equation*}
\frac{1}{2}\sum_{i,j=1}^{d}\frac{\partial ^{2}\mathbf{\mathbf{g}}
	_{k}^{-1}(0;t,Y_{t}\mathbf{)}}{\partial y_{i}\partial y_{j}}dY_{t,i}dY_{t,j},
\text{ }k=1,...,d.
\end{equation*}
With $d\widetilde{Y}_{t}=(d\widetilde{Y}_{t,1},...,d\widetilde{Y}_{t,d})^{T}$
we obtain from (\ref{1}), (\ref{7}) and (\ref{13}):
\begin{equation*}
d\widetilde{Y}_{t}=-\mathbf{g}_{\ast }^{-1}(t;0,\mathbf{\mathbf{g}}
^{-1}(0;t,Y_{t}\mathbf{)})F(t,Y_{t})dt+
\end{equation*}
\begin{equation*}
\mathbf{\mathbf{g}}_{\ast }^{-1}(t;0,\mathbf{\mathbf{g}}^{-1}(0;t,Y_{t}
\mathbf{))}\left[ (F(t,Y_{t})+m(t,Y_{t}))dt+\sigma (t,Y_{t})dW_{t}\right] +
\end{equation*}
\begin{equation*}
\frac{1}{2}\mathbf{\mathbf{g}}_{\ast }^{-1}(t;0,\mathbf{\mathbf{g}}
^{-1}(0;t,Y_{t}\mathbf{))}\sum_{i,j=1}^{d}c_{ij}(t,0,\mathbf{\mathbf{g}}
^{-1}(0;t,Y_{t}\mathbf{)})\sum_{p=1}^{d}[\sigma _{ip}(t,Y_{t})\sigma
_{jp}(t,Y_{t})]dt=
\end{equation*}
\begin{equation*}
\mathbf{\mathbf{g}}_{\ast }^{-1}(t;0,\widetilde{Y}_{t}\mathbf{)}\left\{ m(t,
\mathbf{g(}t;0\mathbf{,}\widetilde{Y}_{t}))+\right. 
\end{equation*}
\begin{equation*}
\left. \frac{1}{2}\sum_{i,j=1}^{d}c_{ij}(t;0,\widetilde{Y}
_{t})\sum_{p=1}^{d}[\sigma _{ip}(t,\mathbf{g(}t;0\mathbf{,}\widetilde{Y}
_{t}))\sigma _{jp}(t,\mathbf{g(}t;0\mathbf{,}\widetilde{Y}_{t}))]\right\} dt+
\end{equation*}
\begin{equation}
\mathbf{\mathbf{g}}_{\ast }^{-1}(t;0,\widetilde{Y}_{t}\mathbf{)}\sigma (t,
\mathbf{g(}t;0\mathbf{,}\widetilde{Y}_{t}))dW_{t}:=\widetilde{m}(t,
\widetilde{Y}_{t})dt+\widetilde{\sigma }(t,\widetilde{Y}_{t})dW_{t},
\label{14}
\end{equation}
where
\begin{equation*}
\widetilde{m}(t,\widetilde{Y}_{t})=\mathbf{\mathbf{g}}_{\ast }^{-1}(t;0,
\widetilde{Y}_{t}\mathbf{)}\left\{ m(t,\mathbf{g(}t;0\mathbf{,}\widetilde{Y}
_{t}))+\right. 
\end{equation*}
\begin{equation*}
\left. \frac{1}{2}\sum_{i,j=1}^{d}c_{ij}(t;0,\widetilde{Y}
_{t})\sum_{p=1}^{d}[\sigma _{ip}(t,\mathbf{g(}t;0\mathbf{,}\widetilde{Y}
_{t}))\sigma _{jp}(t,\mathbf{g(}t;0\mathbf{,}\widetilde{Y}_{t}))]\right\} ,
\end{equation*}
\begin{equation*}
\widetilde{\sigma }(t,\widetilde{Y}_{t}):=\mathbf{\mathbf{g}}_{\ast
}^{-1}(t;0,\widetilde{Y}_{t}\mathbf{)}\sigma (t,\mathbf{g(}t;0\mathbf{,}
\widetilde{Y}_{t})).
\end{equation*}

We obtaine that the process $\widetilde{Y}_{t}$ satisfies the Ito SDE (\ref{14}) with \textit{bounded coefficiens} $\widetilde{m}(t,y)$ and $\widetilde{\sigma }(t,y)$. In the linear case all vectors $c_{jk}=0$ because $\mathbf{g}_{\ast }(t;t_{0},\mathbf{x)}=\Phi (t)$ does not depend on the initial condition $\mathbf{x}$.

	\section{Markov chain}
	
	Suppose that assumptions \textbf{(A1)} -- \textbf{(A3)} hold true.
	
	It follows from the inverse function theorem that $G:R^{d}\rightarrow R^{d}:$ $	x\rightarrow x+hF(t,x)$ is locally a $C^{2}$ -- diffeomorphism for each $t\in \lbrack	0,T]$ and $h$ small enough.
	
	We consider (nonuniform) partition $0=t_{0}^{n}<t_{1}^{n}<...<t_{n}^{n}=T$
	and a Markov chain
	\begin{equation}\label{71}
	X(t_{k+1}^{n})=X(t_{k}^{n})+h_{k}^{n}%
	\{F(t_{k}^{n},X(t_{k}^{n}))+m(t_{k}^{n},X(t_{k}^{n}))\}
	\end{equation}
	\begin{equation*}
	+\sqrt{h_{k}^{n}}\sigma (t_{k}^{n},X(t_{k}^{n}))\cdot \varepsilon
	(t_{k+1}^{n}),
	\end{equation*}
	\begin{equation*}
	X(0)=x\in R^{d},k=0,...,n-1,h_{k}^{n}=t_{k+1}^{n}-t_{k}^{n}. 
	\end{equation*}
	\textbf{(A4)} There exists a constant $C>1$ such that
	\begin{equation*}
	C^{-1}\leq \frac{h_{k}^{n}}{h_{l}^{n}}\leq C,\text{ for }n\geq 1\text{ and }%
	1\leq k,l\leq n.
	\end{equation*}
	\begin{equation*}
	\lim_{n\rightarrow \infty }h_{1}^{n}=0.
	\end{equation*}
	
	Consider the difference equation
	\begin{equation}
	\frac{\widehat{\mathbf{g}}\mathbf{(}t_{k+1}^{n};0\mathbf{,}x)-\widehat{%
			\mathbf{g}}\mathbf{(}t_{k}^{n};0\mathbf{,}x)}{h_{k}^{n}}=F(t_{k}^{n},%
	\widehat{\mathbf{g}}\mathbf{(}t_{k}^{n};0\mathbf{,}x)), \;\;\;\; \widehat{\mathbf{g}}%
	\mathbf{(}0;0\mathbf{,}x)=x,  \label{81}
	\end{equation}
	where $\widehat{\mathbf{g}}\mathbf{(}t_{k}^{n};0\mathbf{,}x)$ s the Euler broken line corresponding to the partition $0=t_{0}^{n}<t_{1}^{n}<...<t_{n}^{n}=T$ starting from $x$ at the moment $0$ and constructed for the solution $\mathbf{g(}t;0\mathbf{,}x)$ of the equaiton $\overset{\cdot }{\mathbf{x}}=F(t,\mathbf{x})\mathbf{,  \;\;g(}0;0\mathbf{,}x)=x.$ If we put $\widehat{X}(t_{k}^{n}):=\widehat{\mathbf{g}}^{-1}\mathbf{%
		(}0;t_{k}^{n}\mathbf{,}X(t_{k}^{n}))$ and replace $x$ by $\widehat{X}%
	(t_{k}^{n})$ in (\ref{81}) we obtain:
	\begin{equation*}
	\frac{\widehat{\mathbf{g}}\mathbf{(}t_{k+1}^{n};0\mathbf{,}\widehat{X}%
		(t_{k}^{n}))-\widehat{\mathbf{g}}(t_{k}^{n};0,\widehat{X}(t_{k}^{n}))}{%
		h_{k}^{n}}=F(t_{k}^{n},\widehat{\mathbf{g}}(t_{k}^{n};0,\widehat{X}%
	(t_{k}^{n}))),
	\end{equation*}%
	\begin{equation*}
	\widehat{X}(t_{k}^{n})=\widehat{\mathbf{g}}(0;0,\widehat{X}(t_{k}^{n}))
	\label{8a}
	\end{equation*}%
	Then (\ref{71}) can be rewritten as
	\begin{equation*}
	\widehat{\mathbf{g}}(t_{k+1}^{n};0,\widehat{X}(t_{k+1}^{n}))=\widehat{%
		\mathbf{g}}\mathbf{(}t_{k+1}^{n};0\mathbf{,}\widehat{X}%
	(t_{k}^{n}))+h_{k}^{n}m(t_{k}^{n},\widehat{\mathbf{g}}(t_{k}^{n};0,\widehat{X%
	}(t_{k}^{n})))+
	\end{equation*}%
	\begin{equation*}
	\sqrt{h_{k}^{n}}\sigma (t_{k}^{n},\widehat{\mathbf{g}}(t_{k}^{n};0,\widehat{X%
	}(t_{k}^{n}))\cdot \varepsilon (t_{k+1}^{n}).  \label{9a}
	\end{equation*}%
	Iterating (\ref{81}) we obtain:
	\begin{equation*}
	\widehat{\mathbf{g}}\mathbf{(}(t_{k+1}^{n};0\mathbf{,}%
	x):=L_{t_{k+1}^{n}}^{n}(x),  \label{10a}
	\end{equation*}%
	where $L_{t}^{n}(x)$ is the Euler broken line starting from $x$ for the
	equation $\overset{\cdot }{\mathbf{x}}=F(t,\mathbf{x})$ corresponding to the partition $0=t_{0}^{n}<t_{1}^{n}<...<t_{n}^{n}=T:$ 
	\begin{equation*}
	L_{t_{1}^{n}}^{n}(x)=x+h_{1}^{n}F(0,x),
	\end{equation*}%
	\begin{equation*}
	...
	\end{equation*}%
	\begin{equation*}
	L_{t_{k+1}^{n}}^{n}(x)=L_{t_{k}^{n}}^{n}(x)+h_{k+1}^{n}F(t_{k}^{n},L_{t_{k}^{n}}^{n}(x))
	\end{equation*}
	
	Iterating the last equation we get:
	\begin{equation}
	\widehat{\mathbf{g}}\mathbf{(}t_{k+1}^{n};0\mathbf{,}x)=x+%
	\sum_{j=0}^{k}h_{j+1}^{n}F(t_{j}^{n},\widehat{\mathbf{g}}\mathbf{(}%
	t_{j}^{n};0\mathbf{,}x)),  \label{11a}
	\end{equation}%
	\begin{equation*}
	\widehat{\mathbf{g}}\mathbf{(}t_{k}^{n};0\mathbf{,}x)=(\widehat{\mathbf{g}}%
	_{1}\mathbf{(}t_{k}^{n};0\mathbf{,}x),...,\widehat{\mathbf{g}}_{d}\mathbf{(}%
	t_{k}^{n};0\mathbf{,}x))^{T}.
	\end{equation*}%
	Define the Jacobian matrix $\widehat{\mathbf{g}}_{\ast}(t_{k}^{n};0,x):=\left\Vert \widehat{z}_{ij}(t_{k}^{n};0,x)\right\Vert , \widehat{z}_{ij}(t_{k}^{n};0,x):=\frac{\partial \widehat{\mathbf{g}}_{i} \mathbf{(}t_{k}^{n};0\mathbf{,}x)}{\partial x_{j}}$. The inverse mapping is
	given by the matrix $\widehat{\mathbf{\mathbf{g}}}_{\ast
	}^{-1}(t_{k}^{n},0,y):=\left\Vert \widehat{z}^{ij}(t_{k}^{n};0,y)\right\Vert ,%
	\widehat{z}^{ij}(t_{k}^{n};0,y):=\frac{\partial \widehat{\mathbf{g}}_{i}^{-1}%
		\mathbf{(}0;t_{k}^{n}\mathbf{,}y)}{\partial y_{j}}$. Differentiating (\ref{11a}) w.r.t. $x$ we obtain:
	\begin{equation}
	\widehat{\mathbf{\mathbf{g}}}_{\ast
	}(t_{k+1}^{n};0,x)=I+\sum_{j=0}^{k}h_{j+1}^{n}F_{\ast }(t_{j}^{n},\widehat{%
	\mathbf{g}}\mathbf{(}t_{j}^{n};0\mathbf{,}x))\widehat{\mathbf{\mathbf{g}}}%
_{\ast }(t_{j}^{n};0,x),  \label{16a}
\end{equation}
\begin{equation}
\left\Vert \widehat{\mathbf{\mathbf{g}}}_{\ast }(t_{k+1}^{n};0,x)\right\Vert
\leq \sqrt{d}+K\sum_{j=0}^{k}h_{j+1}^{n}\left\Vert \widehat{\mathbf{\mathbf{g%
		}}}_{\ast }(t_{j}^{n};0,x)\right\Vert ,  \label{17}
		\end{equation}
		To estimate l.h.s. of (\ref{17}) we apply the following Proposition.
		
		\textbf{Proposition (discrete Gronwall Lemma \cite{KP04})}. Let  $\{y_{n}\}$ and $\{g_{n}\}$ be nonnegative sequences and $c$ a nonnegative constant. If
		\begin{equation*}
		y_{k+1}\leq c+\sum_{j=0}^{k}g_{j}y_{j}, \; k\geq 0,
	\end{equation*}%
	then
	\begin{equation*}
	y_{k+1}\leq c + \prod \limits_{j=0}^{k}(1+g_{j})\leq c\exp \left(
	\sum_{j=0}^{k}g_{j}\right) .
	\end{equation*}%
	We apply this Proposition with $y_{k}=\left\Vert \widehat{\mathbf{\mathbf{g}}%
	}_{\ast }(t_{k}^{n};0,x)\right\Vert ,c=\sqrt{d}%
	,g_{j}=Kh_{j+1}^{n},j=0,...,k.$ We obtain (\ref{17}):
	\begin{equation*}
	\left\Vert \widehat{\mathbf{\mathbf{g}}}_{\ast }(t_{k+1}^{n};0,x)\right\Vert
	\leq \sqrt{d}\exp \left( KT\right) ,k=0,...,n-1.  \label{18}
	\end{equation*}
	By direct substitution we can check that the equation (\ref{16a}) admits an
	explicit solution
	\begin{equation}
	\widehat{\mathbf{\mathbf{g}}}_{\ast
	}(t_{m}^{n};0,x)=\prod\limits_{j=m-1}^{0}[I+h_{j+1}^{n}F_{\ast }(t_{j}^{n},%
	\widehat{\mathbf{g}}\mathbf{(}t_{j}^{n};0\mathbf{,}x))]  \label{19}
	\end{equation}%
	(with the convention $\prod\limits_{j=-1}^{0}[I+h_{j+1}^{n}F_{\ast }(jh,%
	\widehat{\theta }_{jh}(x))]=I$). Indeed, substituting (\ref{19}) in the r.h.s. (\ref{16a}) we obtain a telescopic sum:
	\begin{equation*}
	I+\sum_{j=0}^{k}\left\{ [I+h_{j+1}^{n}F_{\ast }(t_{j}^{n},\widehat{\mathbf{g}%
	}\mathbf{(}t_{j}^{n};0\mathbf{,}x))-I]\prod%
	\limits_{l=j-1}^{0}[I+h_{l+1}^{n}F_{\ast }(t_{l}^{n},\widehat{\mathbf{g}}%
	\mathbf{(}t_{l}^{n};0\mathbf{,}x))]\right\} =
	\end{equation*}
	\begin{equation*}
	I+\sum_{j=0}^{k}\left\{ \prod\limits_{l=j}^{0}[I+h_{l+1}^{n}F_{\ast
	}(t_{l}^{n},\widehat{\mathbf{g}}\mathbf{(}t_{l}^{n};0\mathbf{,}%
	x))]-\prod\limits_{l=j-1}^{0}[I+h_{l+1}^{n}F_{\ast }(t_{l}^{n},\widehat{%
		\mathbf{g}}\mathbf{(}t_{l}^{n};0\mathbf{,}x))]\right\} =
	\end{equation*}%
	\begin{equation*}
	I+\prod\limits_{l=k}^{0}[I+h_{l+1}^{n}F_{\ast }(t_{l}^{n},\widehat{\mathbf{g%
		}}\mathbf{(}t_{l}^{n};0\mathbf{,}x))]-I=\widehat{\mathbf{\mathbf{g}}}_{\ast
	}(t_{k+1}^{n};0,x).
	\end{equation*}%
	If each of the matrices $I+h_{l+1}^{n}F_{\ast }(t_{l}^{n},\widehat{\mathbf{g}%
	}\mathbf{(}t_{l}^{n};0\mathbf{,}x)),$ $l=0,1,...,k$ is invertible then $%
	\widehat{\mathbf{\mathbf{g}}}_{\ast }(t_{k+1}^{n};0,x)$ is also invertible
	and
	\begin{equation*}
	\widehat{\mathbf{\mathbf{g}}}_{\ast
	}^{-1}(t_{k+1}^{n};0,x)=\prod\limits_{j=0}^{k}[I+h_{j+1}^{n}F_{\ast
}(t_{j}^{n},\widehat{\mathbf{g}}\mathbf{(}t_{j}^{n};0\mathbf{,}x))]^{-1}.
\label{20}
\end{equation*}%
By (\textbf{A2}) and \textbf{(A4) } $\left\Vert h_{j+1}^{n}F_{\ast
}(t_{j}^{n},\widehat{\mathbf{g}}\mathbf{(}t_{j}^{n};0\mathbf{,}%
x))\right\Vert \leq h_{j+1}^{n}K<\frac{1}{2}$  for sufficiently large $n$ and, hence, each of the matrices $I+h_{j+1}^{n}F_{\ast }(t_{j}^{n},\widehat{%
	\mathbf{g}}\mathbf{(}t_{j}^{n};0\mathbf{,}x)),$ $j=0,1,...,k$ is invertible.
Hence%
\begin{equation*}
\max_{0\leq m\leq n}\left\Vert \widehat{\mathbf{\mathbf{g}}}_{\ast
}^{-1}(t_{m}^{n};0,x)\right\Vert \leq \max_{0\leq m\leq
n}\prod\limits_{j=0}^{m-1}\left\Vert [I+h_{j+1}^{n}F_{\ast }(t_{j}^{n},%
\widehat{\mathbf{g}}\mathbf{(}t_{j}^{n};0\mathbf{,}x))]^{-1}\right\Vert 
\end{equation*}%
\begin{equation*}
\leq \max_{1\leq m\leq n}\prod\limits_{j=0}^{m-1}\frac{1}{1-\left\Vert
	h_{j+1}^{n}F_{\ast }(t_{j}^{n},\widehat{\mathbf{g}}\mathbf{(}t_{j}^{n};0%
	\mathbf{,}x))\right\Vert }
\end{equation*}%
\begin{equation*}
\leq \max_{1\leq m\leq n}\prod\limits_{j=0}^{m-1}(1+\frac{Kh_{j+1}^{n}}{%
	1-Kh_{j+1}^{n}})\leq (1+2Kh_{j+1}^{n})^{n}\leq C(K,T).  \label{21}
\end{equation*}%
for $n$ large enough. Passing to the limit in (\ref{16}) when $n\rightarrow \infty $ we obtain $\widehat{\mathbf{\mathbf{g}}}_{\ast }(\phi ^{n}(t);0,x)$ $%
\rightarrow \mathbf{g}_{\ast }(t;0,\mathbf{x)},\phi ^{n}(t)=\inf \left\{
t_{i}^{n}:t_{i}^{n}\leq t<t_{i+1}^{n}\right\} ,$ where $\mathbf{g}_{\ast
}(t;0,\mathbf{x)}$ is a solution of
\begin{equation*}
\mathbf{g}_{\ast }(t;0,\mathbf{x)}=I+\int_{0}^{t}F_{\ast }(u,\mathbf{g(}u;0,%
\mathbf{x))g}_{\ast }(u;0,\mathbf{x)}du,
\end{equation*}%
or
\begin{equation*}
\frac{\partial }{\partial t}\mathbf{g}_{\ast }(t;0,\mathbf{x)=}F_{\ast }(t,%
\mathbf{g(}t;0,\mathbf{x))\mathbf{g}_{\ast }(}t;0,\mathbf{\mathbf{x)},%
	\mathbf{g}_{\ast }(}0;0,\mathbf{\mathbf{x)}}=\mathbf{I.}  \label{22}
\end{equation*}%
As we mentioned above $x\rightarrow x+hF(t,x)$ is a $C^{2}$--diffeomorphism for each $t\in \lbrack 0,T]$  and $h$ small enough. It follows that $\widehat{\mathbf{g}}^{-1}\mathbf{(}0;t_{k}^{n}\mathbf{,}y)$ is uniquely defined and
smoothly depends on $y$.Hence the mapping $R^{d}\rightarrow
R^{d}:x\rightarrow \widehat{\mathbf{g}}\mathbf{(}t_{k}^{n};0\mathbf{,}x)$ is a $C^{2}$--diffeomorphism. Differential of this mapping is given by the Jacobian matrix $\widehat{\mathbf{g}}_{\ast }\mathbf{(}t_{k}^{n};0\mathbf{,}%
x):=\left\Vert \widehat{z}_{ij}(t_{k}^{n};0\mathbf{,}x)\right\Vert ,\widehat{%
	z}_{ij}(t_{k}^{n};0\mathbf{,}x)=\frac{\partial \widehat{\mathbf{g}}_{i}%
	\mathbf{(}t_{k}^{n};0\mathbf{,}x)}{\partial x_{j}}.$ Differential of the
inverse mapping is given by the matrix $\widehat{\mathbf{g}}_{\ast }^{-1}%
\mathbf{(}0;t_{k}^{n}\mathbf{,}y):=\left\Vert \widehat{z}^{ij}(0;t_{k}^{n}%
\mathbf{,}y)\right\Vert ,\widehat{z}^{ij}(0;t_{k}^{n}\mathbf{,}y)=\frac{%
	\partial \widehat{\mathbf{g}}_{\ast i}^{-1}\mathbf{(}0;t_{k}^{n}\mathbf{,}y)%
}{\partial y_{j}}$ 
\\
 Define%
\begin{equation*}
\widetilde{X}(t_{k}^{n})=\widehat{\mathbf{g}}^{-1}\mathbf{(}0;t_{k}^{n}%
\mathbf{,}X(t_{k}^{n})).  \label{24}
\end{equation*}%
Then we have:
\begin{equation*}
\widetilde{X}(t_{k+1}^{n})-\widetilde{X}(t_{k}^{n})=\widehat{\mathbf{g}}^{-1}%
\mathbf{(}0;t_{k+1}^{n}\mathbf{,}X(t_{k+1}^{n}))-\widehat{\mathbf{g}}^{-1}%
\mathbf{(}0;t_{k+1}^{n}\mathbf{,}\widehat{\mathbf{g}}%
(h_{k}^{n};t_{k}^{n},X(t_{k}^{n})))=
\end{equation*}%
\begin{equation*}
\left( \int_{0}^{1}\widehat{\mathbf{g}}_{\ast }^{-1}(0;t_{k+1}^{n}\mathbf{,}%
\Psi _{u}(\widetilde{X}(t_{k}^{n}), \varepsilon(t_{k}^{n})))du\right) \left( X(t_{k+1}^{n})-\widehat{%
	\mathbf{g}}(h_{k}^{n};t_{k}^{n},X(t_{k}^{n}))\right) =
\end{equation*}%
\begin{equation*}
h_{k}^{n}\left( \int_{0}^{1}\widehat{\mathbf{g}}_{\ast
}^{-1}(0;t_{k+1}^{n},\Psi _{u}(\widetilde{X}(t_{k}^{n}), \varepsilon(t_{k}^{n})))du\right)
m(t_{k}^{n},\widehat{\mathbf{g}}\mathbf{(}t_{k}^{n};0,\widetilde{X}%
(t_{k}^{n})))
\end{equation*}%
\begin{equation*}
+\sqrt{h_{k}^{n}}\left( \int_{0}^{1}\widehat{\mathbf{g}}_{\ast
}^{-1}(0;t_{k+1}^{n},\Psi _{u}(\widetilde{X}(t_{k}^{n}), \varepsilon(t_{k}^{n})))du\right) \sigma
(t_{k}^{n},\widehat{\mathbf{g}}\mathbf{(}t_{k}^{n};0,\widetilde{X}%
(t_{k}^{n})))\varepsilon (t_{k+1}^{n}),  \label{25}
\end{equation*}%
where 
\begin{equation*}
\Psi _{u}(\widetilde{X}(t_{k}^{n}), \varepsilon(t_{k}^{n})):=\widehat{\mathbf{g}}%
(h_{k}^{n};t_{k}^{n},X(t_{k}^{n}))+u[X(t_{k+1}^{n})-\widehat{\mathbf{g}}%
(h_{k}^{n};t_{k}^{n},X(t_{k}^{n}))]=X(t_{k}^{n})+h_{k}^{n}F(t_{k}^{n},X(t_{k}^{n}))
\end{equation*}%
\begin{equation*}
+uh_{k}^{n}m(h_{k}^{n},X(t_{k}^{n}))+u\sqrt{h_{k}^{n}}\sigma
(h_{k}^{n},X(t_{k}^{n}))\varepsilon (t_{k+1}^{n})=
\end{equation*}

\begin{equation*}
\widehat{\mathbf{g}}\mathbf{(}t_{k}^{n};0,\widetilde{X}%
(t_{k}^{n}))+h_{k}^{n}F(t_{k}^{n},\widehat{\mathbf{g}}\mathbf{(}t_{k}^{n};0,%
\widetilde{X}(t_{k}^{n})))+uh_{k}^{n}m(t_{k}^{n},\widehat{\mathbf{g}}\mathbf{%
	(}t_{k}^{n};0,\widetilde{X}(t_{k}^{n})))+
\end{equation*}%
\begin{equation*}
u\sqrt{h_{k}^{n}}\sigma (t_{k}^{n},\widehat{\mathbf{g}}\mathbf{(}t_{k}^{n};0,%
\widetilde{X}(t_{k}^{n})))\varepsilon (t_{k+1}^{n}).  \label{26}
\end{equation*}

We denote now%
\begin{equation}
\widetilde{m}(t_{k}^{n},\widetilde{X}(t_{k}^{n}), \varepsilon (t_{k+1}^{n}))=\left( \int_{0}^{1}\widehat{\mathbf{g}}%
_{\ast }^{-1}(0;t_{k+1}^{n},\Psi _{u}(\widetilde{X}(t_{k}^{n}), \varepsilon (t_{k+1}^{n})))du\right)
m(t_{k}^{n},\widehat{\mathbf{g}}\mathbf{(}t_{k}^{n};0,\widetilde{X}%
(t_{k}^{n})))  \label{27}
\end{equation}%
\begin{equation}
\widetilde{\sigma }(t_{k}^{n},\widetilde{X}(t_{k}^{n}), \varepsilon (t_{k+1}^{n}))=\left( \int_{0}^{1}\widehat{%
	\mathbf{g}}_{\ast }^{-1}(0;t_{k+1}^{n},\Psi _{u}(\widetilde{X}%
(t_{k}^{n}), \varepsilon (t_{k+1}^{n})))du\right) \sigma (t_{k}^{n},\widehat{\mathbf{g}}\mathbf{(}%
t_{k}^{n};0,\widetilde{X}(t_{k}^{n}))).  \label{28}
\end{equation}%
Then we obtain a Markov chain
\begin{equation}
\widetilde{X}(t_{k+1}^{n})=\widetilde{X}(t_{k}^{n})+h_{k}^{n}\widetilde{m}%
(t_{k}^{n},\widetilde{X}(t_{k}^{n}), \varepsilon (t_{k+1}^{n}))+\sqrt{h_{k}^{n}}\widetilde{\sigma }%
(t_{k}^{n},\widetilde{X}(t_{k}^{n}), \varepsilon (t_{k+1}^{n}))\varepsilon (t_{k+1}^{n})  \label{29}
\end{equation}
with bounded coefficients  $\widetilde{m}(t_{k}^{n},\widetilde{X}%
(t_{k}^{n}),\varepsilon (t_{k+1}^{n}))$ and $\widetilde{\sigma }(t_{k}^{n},\widetilde{X}(t_{k}^{n}),\varepsilon (t_{k+1}^{n})).$%

\textbf{Remark.} Markov Chain (\ref{29}) obtained by application of the procedure of the excluding nonlinear trend is different from the original Markov chain (\ref{71}).The difference is that trend and diffusion now depend on innovation $\varepsilon (t_{k+1}^{n})$. This more general case of Markov chains studied in \cite{BaR16}. These authors considered the following class of Markov chains
$$
X(t_{k+1}^{n})=\Psi_{k}(X(t_{k}^{n}), \frac{\varepsilon_{k+1}}{\sqrt{n}}, \frac{1}{n})
$$
where $\Psi_{k}: \, R^{d}\times R^{N}\times R_{+} \to R^{d}$  are smooth functions, $\varepsilon_{k}, k \in N^{*}$ is a sequence of independent centered random vectors (innovation). The article considered the case $t_{k}^{n}=\frac{k}{n}$. The article considered the convergence $E f(X(t)) \to E f(Y(t))$, where $Y(t)$ is a continuous-time Markov process. The method of the proof is based on the new version of the Malliavin's calculus. Functions $\Psi_{k}$ considered  infinitely differentiable as is usual with this approach. In our case, this corresponds to an infinitely differentiable function $F, m$ and $\sigma$ in (\ref{71}). The main result of \cite{BaR16} obtained on the assumption (1.9) – (1.11) that are not made for the model (\ref{71}) because of an unbounded function $F$.Let $\alpha=\gamma=0$, $\beta=1$ in (1.10) then the condition is not satisfied. For the model (\ref{29}) we can specify a class of models with infinitely differentiable functions $F, m$ and $\sigma$, for which these conditions are satisfied and we can use the results of \cite{BaR16} for this model.

In linear case
\begin{equation*}
F(t,x)=b(t)x,F_{\ast }(t_{j}^{n},\cdot )=b(t_{j}^{n}),
\end{equation*}%
and
\begin{equation}
\widehat{\mathbf{\mathbf{g}}}_{\ast }^{-1}(t_{k+1}^{n};0,\Psi _{u}(%
\widetilde{X}(t_{k}^{n})))=\prod%
\limits_{j=0}^{k}[I+h_{j+1}^{n}b(t_{j}^{n})]^{-1}.  \label{30}
\end{equation}%
If we substitute(\ref{30}) in (\ref{27}) and (\ref{28}) we see that (\ref{29}) is identical to the result of \cite{KMa15} in the case of a uniform grid.

\section{Conclusion}	

This article describes a procedure that allows to reduce the original equation with unbounded trend to an equation with a bounded trend. A similar procedure is given for a Markov chain. This article is a continuation of the work of authors \cite{KMa15} on the modification of the parametrix method. In future it is planned to consider the procedure that allows to reduce the original equation with unbounded diffusion to the equation with the bounded diffusion.

\bibliographystyle{alpha}
\bibliography{bibli}

\end{document}